\documentclass[12pt]{article}

\usepackage{caption}
\usepackage[top=0.75in]{geometry}
\usepackage{amsmath,amsthm,amssymb}
\usepackage{enumerate, xcolor, graphicx}
\usepackage{algorithm}
\usepackage[noend]{algpseudocode}
\usepackage{hyperref}
\usepackage{courier}
\usepackage{lineno}

\newtheorem{thm}{Theorem}

\newtheorem{prop}[thm]{Proposition}
\newtheorem{cor}[thm]{Corollary}

\theoremstyle{definition}

\newtheorem*{rem}{Remark}

\DeclareMathOperator*{\argmax}{argmax}

\begin{document}

\title{The Card Guessing Game: \\ A generating function approach}
\author{Tipaluck Krityakierne$^{1,3}$ \and Thotsaporn Aek Thanatipanonda$^{2,*}$}
\date{%
    \footnotesize{$^1$Department of Mathematics, Faculty of Science, Mahidol University, Bangkok, Thailand\\%
    $^2$Science Division, Mahidol University International College, Nakhon Pathom, Thailand}\\
    $^3$Centre of Excellence in Mathematics, CHE, Bangkok, Thailand\\
    $^*$Correspondence: \href{mailto:thotsaporn@gmail.com}{\texttt{thotsaporn@gmail.com}}\\[2ex]%
}

\maketitle
\thispagestyle{empty}

\begin{abstract}
Consider a card guessing game with complete feedback in which a deck of $n$ cards ordered $1,\dots, n$ is riffle-shuffled once. With the goal to maximize the number of correct guesses, a player guesses cards from the top of the deck one at a time under the optimal strategy until no cards remain. 
We provide an expression for the expected number of correct guesses
with arbitrary number of terms, an accuracy improvement over the results of Liu (2021).
In addition, using generating functions, we give a unified framework for systematically calculating higher-order moments. 
Although the extension of the framework to $k\geq2$ shuffles is not immediately straightforward, we are able to settle a long-standing McGrath’s conjectured optimal strategy described in Bayer and Diaconis (1992) by showing that the optimal guessing strategy for $k=1$ riffle shuffle does not necessarily apply to $k\geq2$ shuffles.\\
\end{abstract}
\hfill{\textit{Dedicated to Doron Zeilberger on the occasion of his 70+1 birthday}}


{\it Doron Zeilberger a.k.a. ``Dr.Z''} is a pioneer of symbolic computation. His Zeilberger’s algorithm, which automates the  closed form calculation for hyper geometric summation (e.g. \cite{Z2,Z1}), was a major breakthrough in the early 90's in the field of computer algebra. Being known for using a computer program as an experimental tool, Dr.Z has demonstrated how the process of mathematical proof could be easier, more accessible, extendable, and perhaps most importantly, enjoyable. We hope that the technique we used in this paper reflects the true spirit of his research.

\begin{enumerate}
\item All theoretical results in this paper, from start to finish, 
were guided by Maple program. 
We first have the slow program empirically count the number of permutations of $n$-card deck after $k$ shuffles. Secondly, we find a recurrence relation for the generating function whose resulting output agrees with that of the slow program for small $n$ and $k.$
Lastly, we obtain the formulas of the $r$th moment about the mean for any $n$-card deck with 1 shuffle.
These formulas are verified to make sure that they agree with the 
results obtained from the recursive program for the generating function.

\item The moment calculus has been applied for a 1-time riffle shuffle. 
We calculate the moments by evaluating the binomial sums
(which also benefited from Zeilberger's algorithm). In the end, we
draw a conclusion about the distribution of our object of interest.

\item Apart from evaluating the expression involving binomial sum, we appeal to the method of undetermined coefficients.
For instance, let us consider the central binomial coefficient $f(n) = \binom{2n}{n}$ which can be expressed through the recurrence relation
\[  (4n+2)f(n)-(n+1)f(n+1) = 0.  \]
Using only the fact that the leading term of $n!$ is 
$\sqrt{2\pi n}\left(\frac{n}{e}\right)^n,$ 
we solve for the series
approximation of $f(n).$
The ansatz for the infinite series of $f(n)$ is 
\[ f(n) = \dfrac{4^n}{\sqrt{\pi n}}
\left(1+\frac{a_1}{n}+\frac{a_2}{n^2} + \dots \right). \]

Substituting this ansatz in the recurrence relation, we end up with system of linear equations for the coefficients $a_1, a_2, \dots$, which can subsequently be used to solve for as many coefficients as we wish.
\end{enumerate}

\section{Introduction}

Many of us dream to visit the famous Monte Carlo  casino, hoping to take advantage of loopholes, winning big, and carry home a lot of money one day. 
While the casino games are designed to make you lose money, Blackjack players can track the profile of cards that have been dealt, and adapt their bet/playing strategies accordingly, so whenever a basic strategy is used by a player, Blackjack's house edge is usually as low as $0.5\%-1\%$. Using a card counting strategy to determine the outcome of the next hand could lead to a player advantage of up to 2\%, which makes Blackjack a high-profile target for card counters since the 1960s.

Although we will not teach you how to make a profit at the casino, in this article, we discuss the optimal strategy to guess the next card in the $n$-card deck with full information after some ``riffle shuffling'' has been performed. That is, a player guesses cards one at a time starting from the top, and the identity of the card guessed is revealed to the player after each guess. The goal is to maximize the total number of correct guesses after all $n$ cards have been drawn. This situation resembles the game of Blackjack as the players can adjust their next guess according to the known information. Relying on a powerful generating function approach, we touch upon the optimal strategy for guessing the next card, and discuss how to compute the expected number, variance, and other higher moments of the number of correct guesses, etc. This enables us to test if the asymptotic normality holds for the distribution of the number of correct guesses.

Recently, Liu \cite{L} investigated the 1-time riffle shuffle case, and obtained the first two terms of the expression for the expected number of correct guesses for a deck of $n$ cards. In Section 2, we extend their results to any higher moments, each with arbitrary number of terms, depending on the desired precision. The information of the moments is then used to summarize the asymptotic distribution. In Section 3, we discuss the $k$-time riffle shuffle case by providing examples that tie together the existing research findings. The difficulties that arise when trying to apply the optimal strategy to $k\geq2$ riffle shuffles lead us to a counterexample for McGrath’s conjectured optimal  strategy  mentioned in  Bayer  and  Diaconis \cite{BD}. Finally, we state open problems, give some applications, and directions for future work.

\subsubsection*{Gilbert-Shannon-Reeds (GSR) Model for Riffle Shuffles}

A riffle shuffle \cite{AZ} is a technique often used in practice when shuffling a deck of cards. 
Consider a deck of $n$ cards labeled consecutively from 1 to $n$ (from top to bottom), facing down on the table. 
Split the deck into two piles in such a way that the probability of cutting the top $t$ cards is $\dfrac{\binom{n}{t}}{2^n}$. Then, interleave the piles back into a single one.
There are $\binom{n}{t}$ ways to interleave the top and bottom piles, and all these possibilities of interleaving are equally likely, so each interleaving has probability $\dfrac{1}{2^n}$ to come up. 

\begin{figure}[h]
    \centering
    \includegraphics[width=0.6\textwidth]{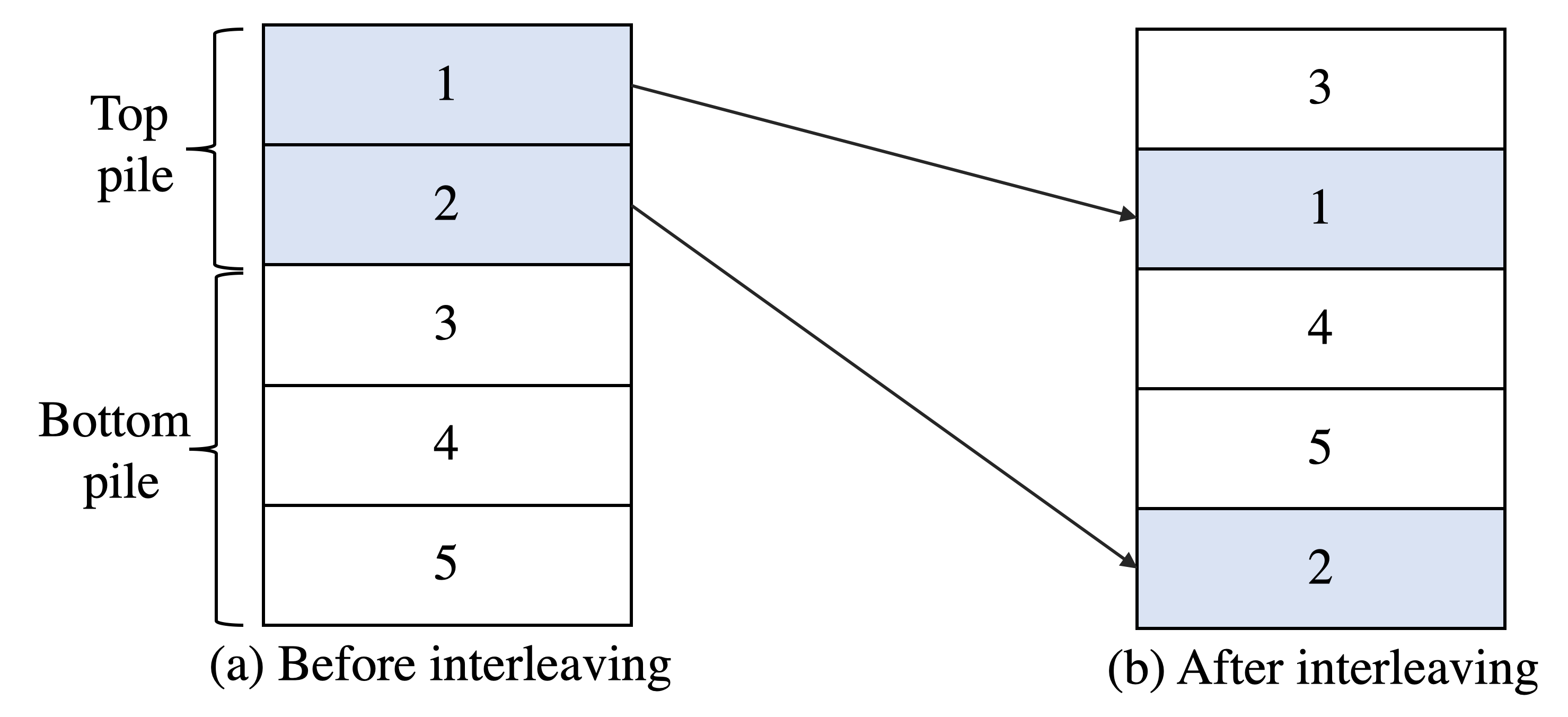}
    \caption{Example of 1-time riffle shuffle of a  deck of 5 cards}
    \label{fig:interleaving}
\end{figure}

For a 1-time shuffle, the operation of interleaving described above gives rise to the identity permutation (having one increasing subsequence) with multiplicity $n+1$, while 
each of the permutations with two increasing subsequences has multiplicity $1.$
And so there are total number of $2^n-n-1$ permutations of two increasing subsequences from the interleaving. An example of 1-time shuffle of a 5-card deck is illustrated in Figure \ref{fig:interleaving}. We can see that there are total of $\binom{5}{2}$ ways to interleave the top and bottom piles.
\section{One-time Riffle Shuffle}

Letting $X_n$ be a random variable representing the number of correct guesses for the $n$-card deck, we define the generating function to keep track of the distribution of $X_n$ in this section. A unified procedure to compute the $r$th factorial moment is subsequently proposed through a recursive formula. Finally, the moments of the distribution are used to summarize the asymptotic distribution of $X_n$ concluding our investigation for the 1-shuffle case.

\subsection{Guessing Strategy and the Number of Correct Guesses}
\subsubsection{Optimal guessing strategy}

The following proposition gives a refresher on the optimal guessing strategy for one-time riffle shuffle \cite{C,L}. It summarizes the optimal strategy for guessing the {\it{first}} and {\it{next card}} {\color{black}(after some number of cards have already been guessed at)}. 

\begin{prop}\label{prop:optimal guessing strategy}
(i) Assume a deck of $n$ cards has been shuffled once. The probability that the first card being $m\in\{1,\dots n\}$ is
\[ p(m) =  \begin{cases}
\dfrac{1}{2}+\dfrac{1}{2^n} & \text{if} \;\  m = 1 \\ \\
\dfrac{\binom{n-1}{m-1}}{2^n}    & \text{if} \;\  m > 1.
\end{cases} \]
Since $p(1)>\frac{1}{2}$, the optimal strategy is to guess card 1 as the first card.

(ii) {\color{black}Assume that the remaining deck} has already been divided into two increasing subsequences of length $a$ and $b$. The probabilities that the next card is from  the subsequences of length $a$ and $b$ are $\dfrac{a}{a+b}$ and $\dfrac{b}{a+b}$, respectively. Hence, the optimal strategy is to guess the first card of the longer subsequence as the next card. 
\end{prop}

{\color{black}While the result $(i)$ of the proposition is precisely the element $M_{11}$ (for $m=1$) and $M_{m1}$ (for $m>1$) of the transition matrix $M$ introduced in Lemma 2.1 of \cite{C}, we give a self-contained alternative proof below. The result $(ii)$ is Proposition 3.1 of \cite{L}.}

\begin{proof}
To prove {\it{(i)}}, we first consider the case when $m=1$. Let $N$ denote the number of permutations whose top card is $1$ (after riffle shuffling once). Depending upon the cutting position $t$,  $N$ can be decomposed into $N=N_0+N_1+\dots N_n$, where $N_t, t=0,1,2,\dots,n$ is the number of permutations whose first card is $1$ given that the cutting position is $t$ before interleaving. Then, $N_0=1$, and $N_t=\binom{n-1}{t-1}$ for $t=1,2,\dots,n$. Thus, the total number of permutations whose first card is 1 (counting all multiplicities) becomes
\[N_0+\sum_{t=1}^n N_t=1+\sum_{t=1}^n \binom{n-1}{t-1}=1+2^{n-1},
\]
and the result for the case when $m=1$ follows.

When $m>1$ {\color{black}(i.e. the cut was made just above card $m$)}, by fixing the first card to be $m$, it is straightforward to see that the number of permutations whose top card is $m$ is $\binom{n-1}{m-1}$, which is the total number of ways to interleave the two subsequences $\{1,2,\dots,m-1\}$ and $\{m+1,\dots n\}$.

{\color{black}The proof of {\it{(ii)}} relies on the fact that the distribution of the resulting interleaving pile of the remaining portion of the deck is essentially uniform (Corollary 2.2 of \cite{L})}, and so the next card comes from the subsequence of length $a$ with probability
\[\frac{{a+b-1\choose a-1}}{{a+b\choose a}}= \frac{a}{a+b}.\]
\qedhere
\end{proof}

\subsubsection{Expected number of correct guesses}
Suppose a player follows the optimal guessing strategy. That is, guess 1 on the first card, and if it turns out that his guess is correct, he continues to guess 2. Now, if the revealed card is not 2, say $m\neq2$, then this implies that {\color{black}the cut was made just above card $m$}, and so the two increasing subsequences at this stage are $\{2,3,\dots,m-1\}$ and $\{m+1,m+2,\dots,n\}$. 
{\color{black}That is, the two increasing subsequences in the remaining
portion of the deck can be determined by the available feedback information.}
His next guess is then either $2$ or $m+1$ depending on the length of the two subsequences. He continues guessing in this way until the deck is completed. Note that if at some point the lengths $a=b$, we assume the player chooses the lower number, e.g. $\min\{2,m+1\}$ in this case. 

\begin{figure}[h]
    \centering
    \includegraphics[width=0.85\textwidth]{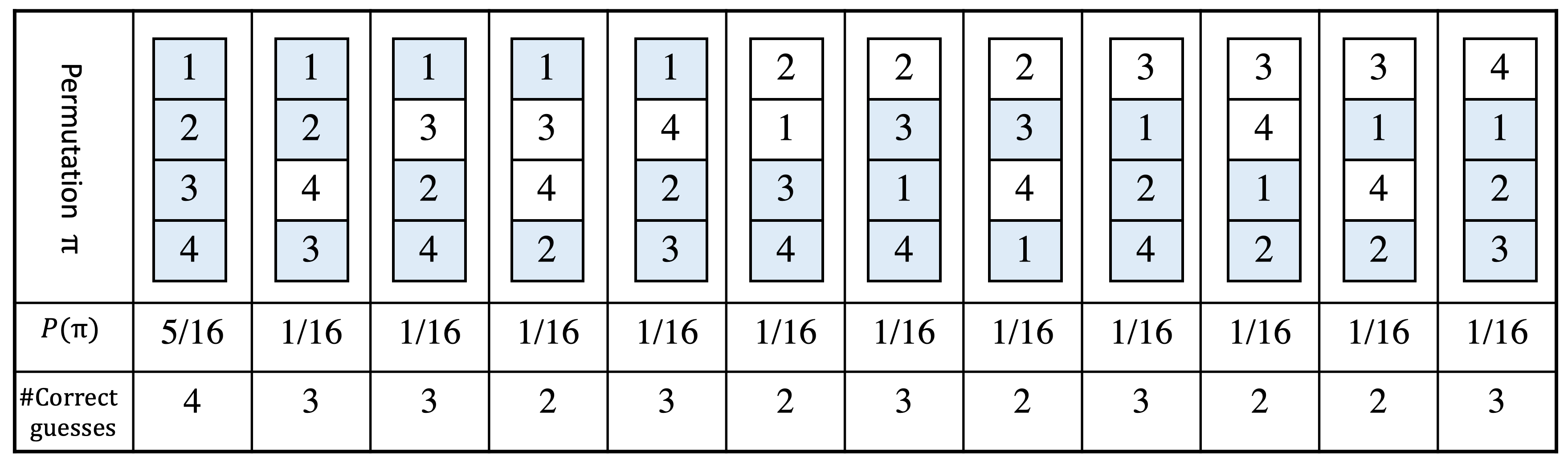}
    \caption{All possible permutations after shuffling a 4-card deck once. The color indicates a correct guess under the optimal strategy.}
    \label{fig:states-4-card-deck}
\end{figure}

Now, using the optimal strategy, what is the expected number of correct guesses? To answer this question, let us examine the distribution of permutations after shuffling a 4-card deck once. Figure \ref{fig:states-4-card-deck} lists all possible permutations along with the probability and the number of correct guesses for a 4-card deck under the optimal strategy. 

Let $X_n$ denote a random variable representing the number of correct guesses (under the optimal strategy) after shuffling an $n$-card deck once. The expected number of correct guesses, $E[X_n]$, can then be computed using the information of the permutation probability distribution. e.g. for a 4-card deck, one can easily obtain
\[
E[X_4]=\sum_{x=2}^4 xP(X_4=x)= \left(2\times\frac{5}{16}\right)+\left(3\times\frac{6}{16}\right)+\left(4\times\frac{5}{16}\right)=3.
\]

\subsection{Generating Functions}
\subsubsection{A recurrence relation for the generating function: motivation}
Let $D_n(q)$ be the (counting) generating function of the number of correct guesses after shuffling a deck of $n$ cards once, i.e. 
\[
    D_n(q)=\sum_{i=0}^{\infty}a_iq^i,    
\]
where $a_i$ denotes the number of 
permutations with $i$ correct guesses. For example, for a 4-card deck, 
\begin{align*} 
D_4(q)&=a_2q^2+a_3q^3+a_4q^4\\  
  &=5q^2+6q^3+5q^4.
\end{align*}
Observe also the connection between $E[X_n]$ and the derivative of $D_n(q)$ w.r.t. $q$:
\[
E[X_n]=\frac{D'_n(q)|_{q=1}}{2^n}.
\]

While our goal is to find a recurrence relation for $D_n(q)$, we will first examine the 4-card deck example, in an attempt to re-express $D_4(q)$ as a function of $D_3(q)$. 

\subsubsection*{Recurrence relation for $D_4(q)$}
\begin{figure}[h]
    \centering
    \includegraphics[width=0.9\textwidth]{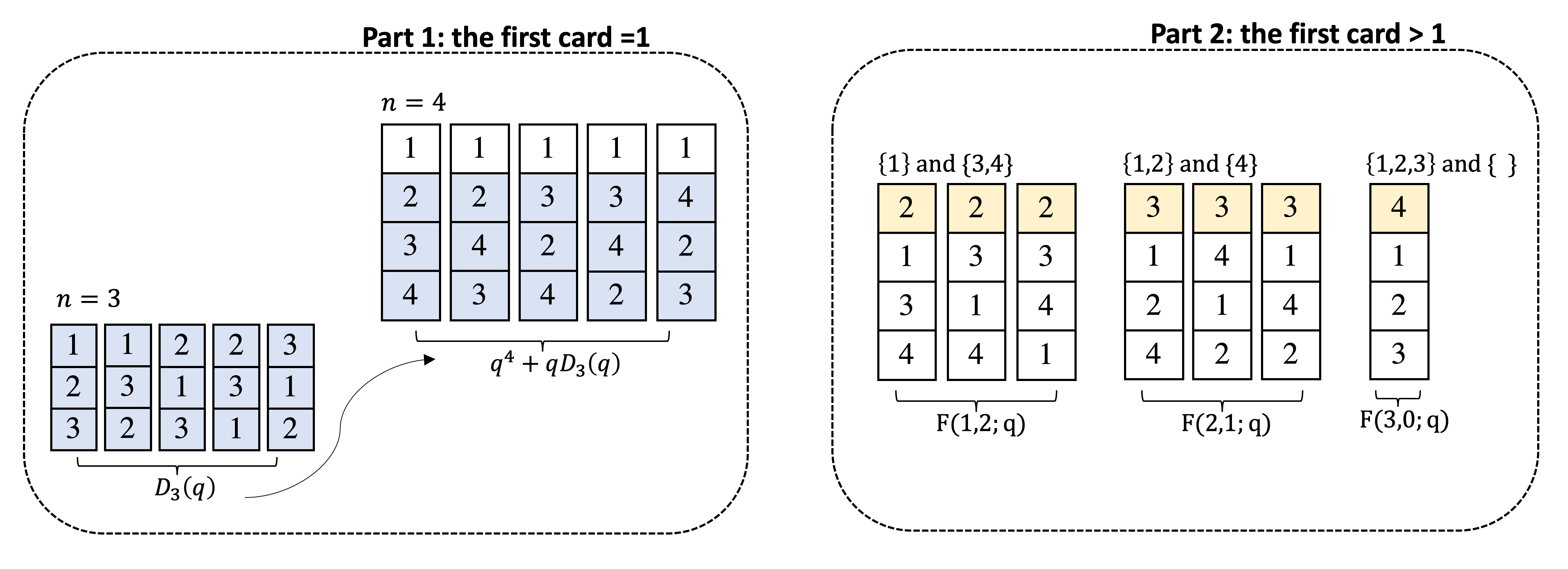}
    \caption{Recurrence structure \scalebox{0.9}{$D_4(q)=\left(q^4+qD_3(q)\right) +F(1,2;q)+F(2,1;q)+F(3,0;q)$}}
    \label{fig:D4recurrence}
\end{figure}

Figure \ref{fig:D4recurrence} illustrates the procedure for finding $D_4(q)$ in a recursive manner. In essence, $D_4(q)$ is decomposed into two parts depending on the identity of the first card.

{\bf{Part 1: the first card being 1.}} Ignoring the first card for a moment, notice how the permutations of the shuffled 3-card deck {\color{black}are} embedded as part of the shuffled 4-card deck. In particular, the correctness of the $i$th guess of the 3-card deck is the same as that of the $(i+1)$th guess of the 4-card deck. Adding back the first card (which is a correct guess under our optimal strategy) yields $qD_3(q)$. Since the number of the identity permutation (having one increasing subsequence) for a shuffled $n$-card deck is $n+1$, it is easy to see why the total contribution to $D_4(q)$ from this part is $q^4+qD_3(q)$. 

{\bf{Part 2: the first card greater than 1.}} As for the other three cases when the first card is greater than 1, guessing the first card under the optimal strategy (i.e. guessing 1) will contribute zero to the total number of correct guesses. However, the identity of the first card {\color{black}lets the player deduce}
as to what the two increasing subsequences are. For example, if the first card is revealed to be 2, the player knows immediately that the two increasing subsequences for the next stage are $\{1\}$ and $\{3,4\}$.  If however 3 is revealed as the first card, then the two increasing subsequences are $\{1,2\}$ and $\{4\}$. Finally, if 4 is revealed as the first card, then the two increasing subsequences are $\{1,2,3\}$ and $\{\,\,\}$, the empty set. In all these cases, the player can guess the next card according to the optimal strategy.  We thus subsequently define a function $F(i,j;q)$ to collect the contribution of the number of correct guesses, where $i$ and $j$ simply represent the lengths of the two increasing subsequences {\color{black} revealed} by the identity of the first card. In particular, using the number of correct guesses given in Figure \ref{fig:states-4-card-deck}, it is obvious that $F(1,2;q)=F(2,1;q)=q^3+2q^2$, and $F(3,0;q)=q^3$. 

Finally, combining the contributions from two parts, the recurrence relation for $D_4(q)$ is expressed as
\[
D_4(q)=\left(q^4+qD_3(q)\right) +F(1,2;q)+F(2,1;q)+F(3,0;q).
\]

\subsubsection{General case}
Now that we have introduced the recurrence relation for $D_4(q)$, we will proceed to define a recursive formula for $D_n(q)$, and go on to find a solution to $F(m,n;q)$ and eventually to $D_n(q)$. 

Along the same lines as $D_4(q)$, the recursive definition of $D_n(q)$ comes from the contributions of two parts:
\begin{equation}
 D_n(q) =\underbrace{qD_{n-1}(q) +q^n}_{\text{the first card $=1$}} +\underbrace{\sum_{i=0}^{n-2} F(n-1-i,i;q)}_{\text{the first card $>1$}},     
\end{equation}
where $D_0(q) = 1$.

By conditioning on the identity of the next card, under the optimal strategy (i.e. guessing the first card of the longer subsequence as the next card), $F(m,n; q)$ can be defined recursively as
\begin{equation} \label{F_recurrence}
  F(m,n; q)  = \underbrace{qF(m-1,n;q)}_{\text{next card from longer subsq.}}+\underbrace{F(m,n-1;q)}_{\text{next card from shorter subsq.}},   
\end{equation}
for $m \geq n$, where $F(m,0; q) = q^m.$ 

{\color{black} 
\begin{rem}
Let us note that under the optimal strategy (which guesses the longer subsequence), the equality $F(m,n; q)=F(n,m; q)$ holds for any $m, n$. In particular, $F(m,n; q)=F(n,m; q)$  whenever $m<n$. Hence, the recurrence relation as defined in \eqref{F_recurrence} indeed covers all cases of $m, n$. 
\end{rem} 
}


The following proposition gives a solution to the recursive function $F(m,n; q)$.
 
\begin{prop} \label{prop:Fmn} 
For $m \geq n,$ 
\begin{equation} \label{F}
  F(m,n;q) = \sum_{i=0}^n 
\left[\binom{m+n}{i}-\binom{m+n}{i-1} \right] q^{m+n-i}.
\end{equation}
\end{prop}

\begin{proof}
When $n=0$, the right hand side of \eqref{F} becomes $q^{m}\binom{m}{0}=q^{m}$ where $\binom{m}{-1}=0$ by convention. {\color{black}We now prove formula \eqref{F} by induction on $m+n$.}

\begin{align*} 
    &F(m,n; q)=qF(m-1,n;q)+F(m,n-1;q)\\
     &=\sum_{i=0}^n \left[\binom{m+n-1}{i}-\binom{m+n-1}{i-1}\right]q^{m+n-i}
    +\sum_{i=0}^{n-1} \left[\binom{m+n-1}{i}-\binom{m+n-1}{i-1}\right]q^{m+n-1-i}\\
    &=q^{m+n}+
     \sum_{i=1}^n \left(\left[\binom{m+n-1}{i}+\binom{m+n-1}{i-1}\right]-\left[\binom{m+n-1}{i-1}+\binom{m+n-1}{i-2}\right]\right)q^{m+n-i}\\
     &=q^{m+n}+
     \sum_{i=1}^n\left[\binom{m+n}{i}-\binom{m+n}{i-1}\right]q^{m+n-i} \textit{ (by Pascal’s triangle rule)} \\
     &=\sum_{i=0}^n\left[\binom{m+n}{i}-\binom{m+n}{i-1}\right]q^{m+n-i}.
\end{align*} 

The first equality is from the recurrence relation \eqref{F_recurrence}
and the second equality is by {\color{black}the induction hypothesis}. 
\end{proof}

\subsection{Moment Calculus}
In this section, we present a derivation of factorial moments in a unified manner, and give expressions for the moments about the mean. The asymptotic distribution of the number of correct guesses will eventually be {\color{black}determined}. For the sake of simplicity, the subscript $n$ of $X$ is omitted in the following analysis.

Let us consider the $r$th factorial moment of $X$ 
\begin{equation} \label{Mo}
 E[X(X-1)\dots(X-r+1)] = \dfrac{ D^{(r)}_n(q)|_{q=1} }{2^n}.
\end{equation}

We first rearrange the recurrence relation of $D_n(q)$
\begin{equation} \label{D}
  D_n(q) = qD_{n-1}(q) + G_n(q), 
\end{equation}
where $G_n(q)$ is defined by
\begin{equation} \label{G}
G_n(q) = q^n+\sum_{i=0}^{n-2} F(n-1-i,i;q). 
\end{equation}

Successively differentiating both sides of \eqref{D} $r$ times ($r \geq 1$) w.r.t. $q$ and substituting in $q=1$, we get 
\begin{equation} \label{Dr}
  D^{(r)}_n(q)|_{q=1} =    D^{(r)}_{n-1}(q)|_{q=1} +
  rD^{(r-1)}_{n-1}(q)|_{q=1} +  G^{(r)}_n(q)|_{q=1}.
\end{equation}

Thus, the $r$th factorial moment can be determined by the following procedure.

\begin{algorithm}
\caption{{\bf Procedure:} Factorial Moment (fixed $r$, formula in $n$)}
\begin{algorithmic}[0]
    \State 1: Compute $G^{(r)}_n(q)|_{q=1}$ by the binomial sum, $n$  symbolic.
    \State 2: Use the method of undetermined coefficient
to calculate $D^{(r)}_{n}(q)|_{q=1}$.
    \State 3: Apply \eqref{Mo} to obtain $E[X(X-1)\dots(X-r+1)]$.
\end{algorithmic}
\end{algorithm}

\subsubsection{Closed-form formula for $G^{(r)}_n(q)|_{q=1}$}

Treating $n$ as a symbolic variable, we consider $G_n(q)$ in \eqref{G} separately for the cases of even and odd $n$.
\[  G_n(q) = \begin{cases}
q^{2k}-q^{2k-1}+2\sum_{j=0}^{k-1} F(2k-1-j, j; q)   &  \text{if} \;\ n=2k, \\ \\
q^{2k+1}-q^{2k}+F(k,k;q)+2\sum_{j=0}^{k-1} F(2k-j, j; q)   &  \text{if} \;\ n=2k+1.
 \end{cases} \]

Applying Proposition \ref{F_recurrence} and interchanging the order of summation,
{\color{black} we obtain
\begin{align*}
\sum_{j=0}^{k-1} F(2k-1-j,j;q) &=
\sum_{j=0}^{k-1} \sum_{i=0}^j
\left[ \binom{2k-1}{i}-\binom{2k-1}{i-1} \right]q^{2k-1-i} \\
&= \sum_{i=0}^{k-1} \sum_{j=i}^{k-1}
\left[ \binom{2k-1}{i}-\binom{2k-1}{i-1} \right]q^{2k-1-i} \\
&= \sum_{i=0}^{k-1} (k-i)\left[ \binom{2k-1}{i}-\binom{2k-1}{i-1} \right]q^{2k-1-i}.
\end{align*}
A similar argument applies to 
$\sum_{j=0}^{k-1} F(2k-j,j;q)$.
Therefore,}
\[  G_n(q) = \begin{cases}
q^{2k}-q^{2k-1}
+2\sum_{i=0}^{k-1} (k-i)
\left[\binom{2k-1}{i}-\binom{2k-1}{i-1} \right]q^{2k-1-i}   & \text{if} \;\ n=2k, \\ \\
 q^{2k+1}-q^{2k}
+2\sum_{i=0}^{k} (k+\frac{1}{2}-i)
\left[\binom{2k}{i}-\binom{2k}{i-1} \right]q^{2k-i} &  \text{if} \;\ n=2k+1.
\end{cases}\]

Hence, the following results follow immediately.

\begin{prop} 
\label{prop:Gr}
For $r\geq1$, the closed-form formula for $G^{(r)}_n(q)|_{q=1}$ can be obtained by evaluating the binomial sums:
\begin{equation} \label{Gr1}  
G^{(r)}_{2k}(q)|_{q=1} = (2k)_r- (2k-1)_r
+2\sum_{i=0}^{k-1} (k-i)
\left[\binom{2k-1}{i}-\binom{2k-1}{i-1} \right] (2k-1-i)_r,  
\end{equation}
\begin{equation} \label{Gr2}  
 G^{(r)}_{2k+1}(q)|_{q=1} = (2k+1)_r-(2k)_r
+2\sum_{i=0}^{k} (k+\frac{1}{2}-i)
\left[\binom{2k}{i}-\binom{2k}{i-1} \right](2k-i)_r, 
\end{equation}
where $(a)_r$ is the falling factorial, i.e. $ (a)_r = a(a-1)(a-2) \dots (a-r+1).$
\end{prop}

We can evaluate the above binomial summations in closed form and obtain a simpler expression for $G^{(r)}_n(q)|_{q=1}$. This can be done either by hand (which could be overwhelming) or by a computer program. Here, we give an example for zero- and first- order moments.
\subsubsection{Zero- and first-order moments}

The zeroth moment corresponds to the number of ways to do one-time riffle shuffle. {\color{black} Using \eqref{F}, we deduce} that $F(m,n;1) = \binom{m+n}{n}$, and a straightforward computation from \eqref{G} gives
\[ 
G_n(1) = 1+\sum_{i=0}^{n-2} \binom{n-1}{i} = 2^{n-1}.
\]
From \eqref{D}, we get the recurrence
\[  
D_n(1) = D_{n-1}(1)+2^{n-1} 
\]
whose solution is 
\[  
D_n(1) = 2^{n}. 
\]
To obtain the first moment, we follow the three-step procedure. When $r = 1$, \eqref{Dr} becomes
\begin{equation} \label{Dr11}  
D'_n(q)|_{q=1} = D'_{n-1}(q)|_{q=1}
+D_{n-1}(1)+ G'_n(q)|_{q=1},   
\end{equation}
where $D_{n-1}(1)=2^{n-1}$ by the zeroth moment.

Maple program {\bf{CloseGr(1,k)}} in Appendix B is used to evaluate the binomial summation
in \eqref{Gr1} and \eqref{Gr2} to get
\begin{align*} 
G'_{2k}(q)|_{q=1} &=  \dfrac{k-1}{2}4^{k}
+ k\binom{2k}{k}+1,  \\
G'_{2k+1}(q)|_{q=1} &= \dfrac{2k-1}{2}4^{k}
+\dfrac{4k+1}{2}\binom{2k}{k}+1.  
\end{align*}

Substituting these two equations back to \eqref{Dr11} separately for odd and even $n$, $D'_{2k+1}(q)|_{q=1}$ and $D'_{2k}(q)|_{q=1}$ can be determined.

In particular, for $n = 2k+1$, 
\[ D'_{2k+1}(q)|_{q=1} = D'_{2k-1}(q)|_{q=1} + \frac{k}{2}4^{k}+ k\binom{2k}{k}+1
 + \frac{2k+1}{2}4^{k}+\dfrac{4k+1}{2}\binom{2k}{k}+1.  \]
 
Simplifying the equation and writing expression in terms of $n$,
\begin{equation}\label{eq:Dr11_odd}
D'_{n}(q)|_{q=1} = D'_{n-2}(q)|_{q=1} 
 + \frac{3n-1}{8}2^{n}+\dfrac{3n-2}{2}\binom{n-1}{(n-1)/2} + 2.     
\end{equation}

Subsequently, $D'_{n}(q)|_{q=1}$ is solved by the method of undetermined coefficients, i.e.
for the term $\dfrac{3n-1}{8}2^{n}$, we assume the solution to be in the form
\begin{equation}  \tag{ansatz1}
(an+b)2^n.
\end{equation}

It is known however that a closed form formula does not exist for $\sum_{k=0}^n P(k)\binom{2k}{k}$, where $P(k)$ is a fixed polynomial in $k$, i.e. not Gosper-summable. Hence, for the term $\dfrac{3n-2}{2}\binom{n-1}{(n-1)/2}$,  
we assume a solution of the form 
\begin{equation}  \tag{ansatz2}
\left(a_0\sqrt{n}+\dfrac{a_1}{\sqrt{n}}+\dfrac{a_2}{n^{3/2}}+\dots\right)2^n.
\end{equation}
Infinite series of the central binomial is given in Appendix A for reference.

The constant term in  \eqref{eq:Dr11_odd} 
{\color{black} and the general solution $C_1+C_2(-1)^n$ of 
$D'_{n}(q)|_{q=1} = D'_{n-2}(q)|_{q=1}$ are}
asymptotically negligible. We thus solve the two parts separately, and then combine their results. 

The final solution (from Maple program {\bf{ForDr(1,1,5,n)}} in Appendix B) is
\[  D'_n(q)|_{q=1} = (n-1)2^{n-1}+ 2^n\sqrt{\frac{2n}{\pi}}\left( 1 - \frac{3}{4n} 
- \frac{53}{96n^2}- \frac{443}{384n^3}- \frac{75949}{18432n^4}
- \frac{4621519}{221184n^5} - \dots \right).\]

Now as $E[X] = \dfrac{D'_n(q)|_{q=1} }{2^n}$, the expectation for the case when $n=2k+1$ is given by
\[ E[X] =   \dfrac{n}{2}+\sqrt{\frac{2n}{\pi}}-\dfrac{1}{2}
-\sqrt{\frac{2}{\pi n}}\left( \frac{3}{4} 
+ \frac{53}{96n} + \frac{443}{384n^2}
+ \frac{75949}{18432n^3}
+ \frac{4621519}{221184n^4} + \dots \right) .\]

In a similar manner, we obtain the expected value for the case when $n=2k$:
\[ E[X] =   \dfrac{n}{2}+\sqrt{\frac{2n}{\pi}}-\dfrac{1}{2}
-\sqrt{\frac{2}{\pi n}}\left( \frac{3}{4} 
+ \frac{49}{96n} + \frac{439}{384n^2}
+ \frac{76709}{18432n^3}
+ \frac{4628519}{221184n^4} + \dots \right) .\]

Note that these results improve the following main result of \cite{L}
\[
E[X]= \dfrac{n}{2}+\sqrt{\dfrac{2n}{\pi}}+O(1).
\]

\newpage
{\color{black}
\textbf{The calculation behind the particular solution for $\displaystyle\dfrac{3n-2}{2}\binom{n-1}{(n-1)/2}$} 

We use this opportunity to show, as an illustration, the calculation of the first few terms of the particular solution of the term $\displaystyle\dfrac{3n-2}{2}\binom{n-1}{(n-1)/2}$. While most textbooks on elementary differential equations discuss the method of undetermined coefficients (e.g. \cite[Chapter 3.5]{BDM}), it is worth the extra effort to discuss this in more detail here, especially the calculations for those terms involving the one-half power.

We start by reminding the reader the Taylor series expansion at $n=\infty$ for the terms such as $\sqrt{n-2}$ and $\dfrac{1}{\sqrt{n-2}}$:
\[ \sqrt{n-2} = \sqrt{n}-\dfrac{1}{\sqrt{n}}-\dfrac{1}{2n^{3/2}}+\dots , \]
\[ \dfrac{1}{\sqrt{n-2}} = \dfrac{1}{\sqrt{n}}+  
\dfrac{1}{n^{3/2}}+
\dfrac{3}{2n^{5/2}}+ \dots .
\]

Now, let us solve the relation
\begin{equation}  \label{TT} D'_{n}(q)|_{q=1} =D'_{n-2}(q)|_{q=1} +\dfrac{3n-2}{2}\binom{n-1}{(n-1)/2}. 
\end{equation}

We first calculate the expansion of the binomial term (using Appendix A):
\begin{small}
\begin{align*}  
&\dfrac{3n-2}{2}\binom{n-1}{(n-1)/2}\\ 
&= \dfrac{3n-2}{2}\cdot\dfrac{1}{\sqrt{\pi}}\left( \dfrac{1}{\sqrt{(n-1)/2}}-\dfrac{1}{8((n-1)/2)^{3/2}}+\dfrac{1}{128((n-1)/2)^{5/2}}+\dots \right)4^{(n-1)/2} \\
&= \dfrac{3n-2}{2}\cdot\dfrac{1}{\sqrt{\pi}}\left(
\sqrt{2}\left( \dfrac{1}{\sqrt{n}}+\dfrac{1}{2n^{3/2}}+\dfrac{3}{8n^{5/2}}+\dots  \right)
-\dfrac{2^{3/2}}{8}\left( \dfrac{1}{n^{3/2}}+\dfrac{3}{2n^{5/2}}+\dfrac{15}{8n^{7/2}}+\dots  \right)+\dots \right)2^{n-1}, 
\end{align*}
\end{small}
where the last equality follows from the Taylor series of $\displaystyle\dfrac{1}{\sqrt{n-1}}$ and $\displaystyle\dfrac{1}{(n-1)^{3/2}}$.

We see that it makes sense to search for a solution
in the form of (ansatz2):
\begin{equation}
\label{eq:Dprime_ans2}
    D'_{n}(q)|_{q=1} =  \left(a_0\sqrt{n}+\dfrac{a_1}{\sqrt{n}}+\dfrac{a_2}{n^{3/2}}+\dots\right)2^n.
\end{equation}
It follows that
\begin{small}
\begin{align*} D'_{n-2}(q)|_{q=1} &=  \left(a_0\sqrt{n-2}+\dfrac{a_1}{\sqrt{n-2}}+\dfrac{a_2}{(n-2)^{3/2}}+\dots\right)2^{n-2} \\
&= a_0\left( \sqrt{n}-\dfrac{1}{\sqrt{n}}-\dfrac{1}{2n^{3/2}}+\dots \right)2^{n-2} +a_1\left( \dfrac{1}{\sqrt{n}}+  
\dfrac{1}{n^{3/2}}+
\dfrac{3}{2n^{5/2}}+ \dots \right)2^{n-2} +\dots .
\end{align*}
\end{small}

We proceed to compare the coefficients of the same terms and solve for $a_0$ and $a_1$.
Substituting the expressions of $\displaystyle D'_{n-2}(q)|_{q=1}$ and $\displaystyle\dfrac{3n-2}{2}\binom{n-1}{(n-1)/2}$ into \eqref{TT}, and comparing the coefficient of $\sqrt{n}\cdot 2^n$ appearing in \eqref{TT} and \eqref{eq:Dprime_ans2}, we obtain
\[ a_0 = \dfrac{a_0}{4}+\dfrac{3\sqrt{2}}{4\sqrt{\pi}}.\]
Hence, $a_0 = \sqrt{\dfrac{2}{\pi}}$.

Similarly, by comparing the coefficient of $\dfrac{1}{\sqrt{n}}\cdot 2^n$ in \eqref{TT} and \eqref{eq:Dprime_ans2}, we obtain
\[ a_1 = \dfrac{a_1-a_0}{4}-\dfrac{\sqrt{2}}{2\sqrt{\pi}}+\dfrac{3}{4\sqrt{\pi}}\left(  \dfrac{\sqrt{2}}{2}-\dfrac{1}{2\sqrt{2}}\right),\]
and so $a_1 = -\dfrac{3}{4} \sqrt{\dfrac{2}{\pi}}$.

As the number of terms increases, the calculations inevitably become more tedious. Therefore, we use the symbolic computation program to find other coefficients in the expansions.
}
\subsubsection{Higher factorial moments}
Along the same line as that of the first moment, higher factorial moments can be obtained following the three-step procedure. We give expressions to $E[X(X-1)\dots(X-r+1)]$ for $r=2,3$ below. 


For odd $n$:
\[ E[X(X-1)] =   \dfrac{n^2}{4}+\sqrt{\frac{2}{\pi}}n^{3/2}-\dfrac{n}{4}
-\sqrt{\frac{2n}{\pi}}\left( \frac{11}{4} 
+ \frac{5}{96n} + \frac{1753}{1152n^2}
+ \frac{13733}{2048n^3}+ \dots \right) .\]

For even $n$:
\[ E[X(X-1)] =   \dfrac{n^2}{4}+\sqrt{\frac{2}{\pi}}n^{3/2}-\dfrac{n}{4}
-\sqrt{\frac{2n}{\pi}}\left( \frac{11}{4} 
+ \frac{1}{96n} + \frac{1901}{1152n^2}
+ \frac{13917}{2048n^3}+ \dots \right) .\]

For odd $n$:
\[ E[X(X-1)(X-2)] =   \dfrac{n^3}{8}+\frac{3}{4}\sqrt{\frac{2}{\pi}}n^{5/2}
-\frac{65}{16}\sqrt{\frac{2}{\pi}}n^{3/2} -\dfrac{13n}{8}
+\sqrt{\frac{2n}{\pi}}\left( 
 \frac{763}{128} - \frac{2681}{1536n}
- \frac{443239}{73728n^2}+ \dots \right) .\]

For even $n$:
\[ E[X(X-1)(X-2)] =   \dfrac{n^3}{8}+\frac{3}{4}\sqrt{\frac{2}{\pi}}n^{5/2}
-\frac{65}{16}\sqrt{\frac{2}{\pi}}n^{3/2} -\dfrac{13n}{8}
+\sqrt{\frac{2n}{\pi}}\left( 
 \frac{767}{128} - \frac{3085}{1536n}
- \frac{413903}{73728n^2}+ \dots \right) .\]
\subsubsection{Moments about the mean}

To complement the expression of the $r$th factorial moment, we find explicit expressions for the moments about the mean. 
The $r$th raw moment (which later leads to the moment about the mean) can be computed from the
factorial moments through the relation:
\[  E[X^r] = \sum_{i=0}^r {\displaystyle  \left\{{r \atop i}\right\}} E[(X)_i] , \]
where ${\displaystyle  \left\{{r \atop i}\right\}}$ is the Stirling number of the second kind.

For the sake of brevity, we give examples for the case when $n$ is even, and $r=2,3$.

When $n$ is even,
\[ E[(X-\mu)^2] =   \left( \frac{3}{4}-\frac{2}{\pi} \right)n
-\frac{3}{4}+\frac{3}{\pi} -\sqrt{\frac{2}{\pi n}} 
+\frac{11}{12\pi n}+\dots.\]
\begin{align*}
 & E[(X-\mu)^3] = \\
 & \sqrt{\frac{2}{\pi}} 
\left(   \left( \frac{4}{\pi}-\frac{5}{4} \right)n^{3/2}
+ \left( \frac{43}{16}-\frac{9}{\pi} \right)n^{1/2} 
-\frac{3\sqrt{2\pi}}{4}+3\sqrt{\frac{2}{\pi}} 
-\left(  \frac{241}{128}
-\frac{5}{8\pi} \right)\frac{1}{\sqrt{n}}+\dots   \right). 
\end{align*}

The expressions up to order 12 for both even and odd $n$ along with the Maple program used to generate $E[(X-\mu)^r]$ and $E[X(X-1)\dots(X-r+1)]$ are provided at \small{\url{https://thotsaporn.com/Card1.html}}. A brief description of the Maple programs used to evaluate results of this paper is given in Appendix B.

\subsection{Asymptotic Distribution of the Number of Correct Guesses}

The skewness coefficient is given by $\frac{m_3}{m_2^{3/2}}$, where $m_r := E[(X-\mu)^r]$ is the $r$th moment about the mean. 
Using the expression for $m_r$ derived in the previous section, we see that the skewness of $X_n$ does not tend to zero (and hence not symmetric). 
Therefore, the number of correct guesses is not asymptotically normally distributed \cite{Z4}. 

\begin{figure}[h!]
    \centering
    \includegraphics[width=0.75\textwidth]{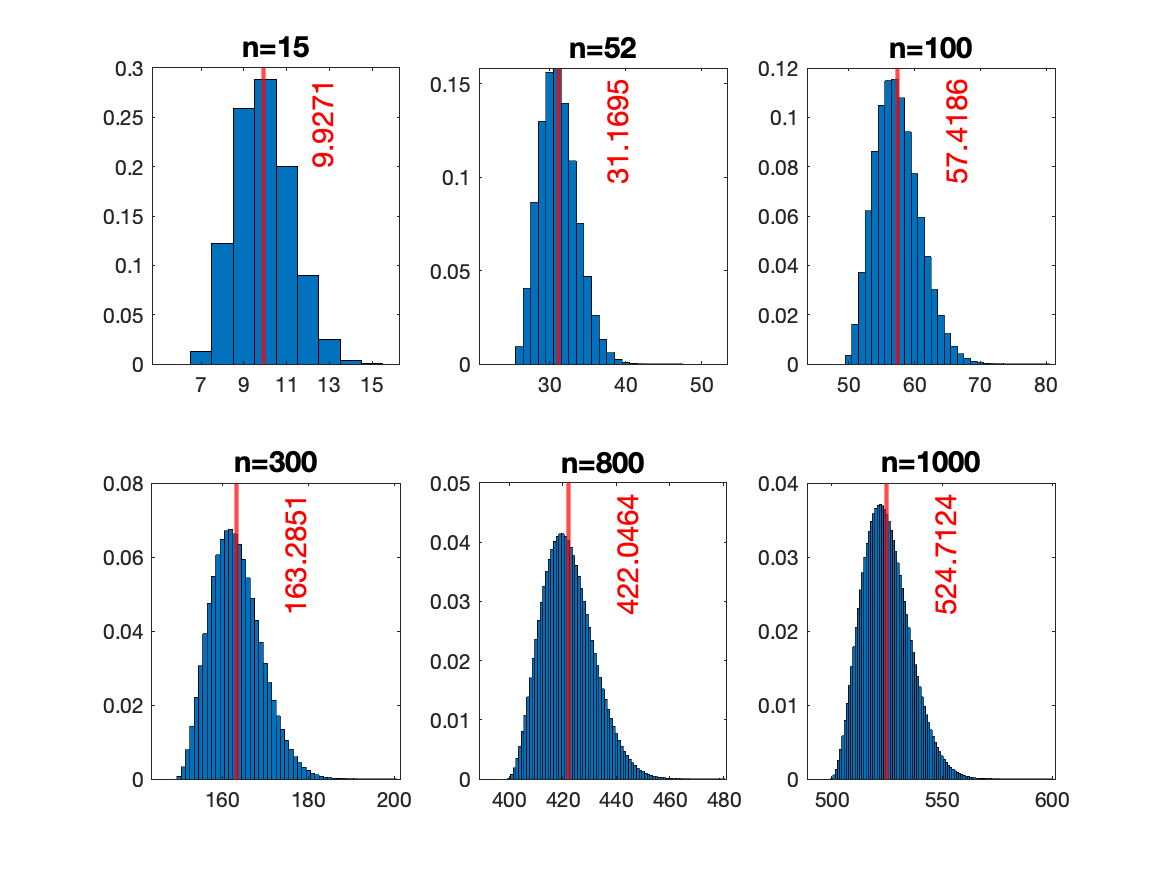}
    \caption{Probability histograms of the number of correct guesses when $n$ varies. The red vertical line indicates the corresponding expected value $E[X_n]$.}
    \label{fig:histogram}
\end{figure}

Histograms illustrating the distribution of $X_n$ for various $n$ are shown in Figure \ref{fig:histogram}. We see clearly that the histogram has right-skewed distribution inconsistent with asymptotic normality.

\newpage
\section{Generalization to $k$ riffle shuffles: is it possible?}

\textbf{Combinatorial interpretation}

We see earlier that shuffling a deck of $n$ card once, the identity permutation (one increasing subsequence) has a multiplicity of $n+1$, while a permutation with two increasing subsequences has a multiplicity $1$. In fact, the multiplicity of a permutation $\pi$ after $k$ shuffles is closely related to the number of increasing subsequences in $\pi$. We state without proof the following result established in \cite{BD}.

\begin{prop} [Theorem 1 of \cite{BD}]  \label{Eu}
Let the permutation $\pi$ of length $n$ 
be the resulting deck after $k$ shuffles.
Let $m$ be the number of increasing 
subsequences of $\pi$, \, $1 \leq m \leq 2^k$. 
The multiplicity of $\pi$ after $k$ shuffles is
\[ \binom{n+2^k-m}{n}. \] 
\end{prop}

{\bf{Example 1}} (Multiplicity of $\pi$){\bf.} Shuffling a deck of 3 cards twice ($n=3$, $k=2$),  \vspace{-0.5em}
\begin{itemize}
    \item the multiplicity of $[1,2,3]$ is $ \binom{3+4-1}{3} = 20;$ \vspace{-0.5em}
    \item the multiplicity of each of $[1,3,2], [2,1,3], [2,3,1], [3,1,2]$ is  $ \binom{3+4-2}{3} = 10;$\vspace{-0.5em}
    \item the multiplicity of $[3,2,1]$  is $\binom{3+4-3}{3} = 4.$ 
\end{itemize}

Thus, the total number of permutations (including the multiplicities) is $(1\times20)+(4\times10)+(1\times4)=2^6$. It is easy to see that in general there are $2^{kn}$ permutations produced by $k$ shuffles of an $n$-card deck. 

We now digress temporarily to present some results concerning $\pi^{-1}$, the inverse permutation, through the following illustration.

{\bf{Example 2}}  (Inverse permutation, $\pi^{-1}$){\bf.} Consider $\pi=[5, 3, 4, 6, 1, 7, 2, 8, 9, 10]$ which has three increasing 
subsequences, namely, $[5,6,7,8,9,10]$, $[3,4]$ and $[1,2]$. Now, since $\pi^{-1} = [5,7,2,3,1,4,6,8,9,10],$ it leads to two {\color{black}descents}, namely, $[7,2], [3,1]$
and three permutation runs, $[\underbrace{5,7}_{\text{run 1}}|\underbrace{2,3}_{\text{run 2}}|\underbrace{1,4,6,8,9,10}_{\text{run 3}}]$. Essentially, if $\pi$ is a permutation after $k$ shuffles with $m$ increasing subsequences, then $\pi^{-1}$ is the permutation with $m-1$ {\color{black}descents} and $m$ runs.

Through the concept of $m$ runs of the inverse riffle shuffles, one can give a combinatorial proof (for the special case, $x=2^k$) of Worpitzky's identity {\color{black}\cite{W}}:
\[ x^n = \sum_{m=1}^{x}  A(n,m-1)\binom{n+x-m}{n}, \]
where $A(n,r)$ is the number of permutations of length $n$ with $r$ descent positions. 
This special case $x=2^k$ for some $k$ gives a formula to the total number of permutations produced by $k$ shuffles of an $n$-card deck.

{ \color{black}
\textbf{Optimal guessing strategy for the $k$ riffle shuffles: a computational approach}

While Proposition 1 gives a clear strategy of how to guess the cards for one shuffle, things are not so clear for $k$ shuffles, $k \geq 2$. 
To collect the data, we apply the program {\bf Best Guess} to {\color{black}sequentially output the best guess at each stage (after each card is revealed)} for an $n$-card deck after $k$ riffle shuffles. 
}

{\color{black}
\begin{algorithm}
\caption{\color{black}{\bf Program:} Best Guess (fixed $n, k$)}
\color{black}
\begin{algorithmic}[0]
    \State \textsc{Goal:} Find the most likely number to show up next
    \State \textsc{Input:} The list of the already revealed cards, $S$
    \State \textsc{Compute:} The number of ways that the number $i$ will show up as the next card
\begin{equation}
\label{eq:Ni}
    N_i = \sum_{\pi \in \Pi_n\left[S,i\right]} \binom{n+2^k-m_{\pi}}{n},
\end{equation} 
where $\Pi_n\left[S,i\right]$ is the set containing all permutations of length $n$ that start with elements in $S$, appended by number $i$. Here, $m_{\pi}$ is the number of increasing subsequences of $\pi$.
\State \textsc{Output:} The best guess $$\displaystyle i^*:=\argmax_i  N_i$$
\end{algorithmic}
\end{algorithm}
}



{\color{black}
The binomial in Equation \eqref{eq:Ni} represents the multiplicity of $\pi$ after $k$ shuffles (Proposition \ref{Eu}). The number $i^*$ will be the optimal guess for the next card.

The challenge of writing the program is to quickly generate $\pi$ (starting with $S$ and $i$) along with counting the number of increasing subsequences of it. 
For this purpose, two versions of the above program were implemented to validate the results. The first program, {\bf SlowBG(n,k,S)}, counts $\pi$ in a rather straightforward way, i.e. list all the permutations in $\Pi_n\left[S,i\right]$ that begin with $S$, prior to starting calculation. The faster program {\bf FastBG(n,k,S)}, on the other hand, generates a permutation which begins with elements $S$ and $i$ in a recursive manner (i.e. appending a number to the permutation one at a time) while keeping track of the number of increasing subsequences at each stage. 
}

Having to calculate the optimal guess at each stage {\color{black}seems to make it
impossible} to set up a convenient recurrence formula for the generating functions of the number of correct guesses (as in the case of 1-time shuffle).  {\color{black}We leave this as an open question.}

\textbf{Open Question:}
{\it Find the optimal strategy for guessing the next card of the $k$ riffle shuffles, $k \geq 2$.}

\textbf{A counterexample}

{\color{black} In the process of searching for the optimal strategy via {\bf SlowBG(n,k,S)/}{\bf FastBG(n,k,S)}, we were able to settle a conjectured optimal strategy due to McGrath described three decades ago in Bayer and Diaconis \cite[Section 5.1]{BD}, which assumed the optimal strategy for a 1-time shuffle applies also for $k\geq2$. The conjecture later appeared
as the eighth problem on a list of ten unsolved problems in Diaconis \cite[Section 5]{D}. We give a concrete
counterexample to disprove McGrath’s conjectured optimal strategy in the following corollary.}

{\color{black} 
\begin{cor}
The optimal guessing strategy for $k$-time riffle shuffles needs not be the same as the optimal strategy for a 1-time riffle shuffle (i.e. guessing the next card from the longer subsequence) as in Proposition \ref{prop:optimal guessing strategy}.
\end{cor}
}

\begin{proof}
{\color{black} We shall give a counterexample for this purpose.} Consider a deck of 10 cards that has been shuffled twice ($n=10$, $k=2$). Now, suppose it is revealed that the first card is 5. Following the conjectured optimal strategy, one would guess 6 as the next card because the remaining subsequences are \{1,2,3,4\}
and \{6,7,8,9,10\} with the latter being a longer subsequence. 
However, careful calculations show that the probability that 1 is the next card is
$\dfrac{31752}{105336}$, while  
the probability that 6 is the next card is
$\dfrac{31570}{105336}.$ {\color{black} The command {\bf SlowBG(10,2,[5])} or {\bf FastBG(10,2,[5])}  outputs these two probabilities (along with the probabilities that card number $i$ will be the next card after 5, for $i=1,\dots, 10$).}
\end{proof}

The above example with $n=10$ and $k=2$, in fact, is the smallest counterexample, i.e. the smallest $n$ for 2 shuffles. {\color{black}As the total number of permutations of this example is $2^{2\times10}$ and the difference between the two probabilities is as small as $\frac{182}{105336}$, it is obvious that finding such a counterexample is a non-trivial task.} 

\newpage
\section{Last But Not Least} 

\textbf{Indicator for the randomness of the deck} 

Perhaps one of the most popular questions that people ask around this topic is ``Is my deck random enough?''. One way to answer this question is to examine the randomness of the shuffled deck through the expected number of correct guesses and compare it to the expected number of correct guesses for a completely random deck. 

For example, let us consider a deck of two cards.  Let $e_k$ denote the expected number of correct guesses as a function of shuffle number $k$. Then, $e_k,\;\ k =0,1,2,\dots$, {\color{black}turn out to be}
\[  2, \frac{7}{4},  \frac{13}{8},  \frac{25}{16},  
\frac{49}{32}, \frac{97}{64}, {\color{black}\frac{193}{128}, \frac{385}{256}, \frac{769}{512},} \dots. \]

On the other hand, the expected number of correct guesses for a randomly shuffled deck is
\[ E[Y] = \frac{1}{2}+1 = \frac{3}{2}. \]

Thus, using $E[Y]$ as a target value, if the threshold value is chosen, say, as 5\% on the relative residual, i.e. $\left|\frac{e_k-1.5}{1.5}\right| <0.05$ in our case, then we need to shuffle the deck 3 times.

Now, for the deck of three cards, we have that $e_k,\;\ k =0,1,2,\dots$ , {\color{black}are}
\[
     3, \frac{19}{8},  \frac{67}{32},  \frac{251}{128},  
\frac{971}{512}, \frac{3819}{2048}, {\color{black}\frac{15147}{8192}, \frac{60331}{32768}, \frac{240811}{131072},} \dots, 
\]

while \[ E[Y] = \frac{1}{3}+\frac{1}{2}+1 = \frac{11}{6}. \]

With the 5\% threshold from the target value, we will need to shuffle the deck 4 times. This idea was discussed before in Section 5.1, a \textit{different distance}, of \cite{BD}.

{\color{black} \textbf{The calculation behind the sequence $e_k$ above} 

Consider a deck of $n$ cards.  Let $\pi_s$ be a permutation of length $n$, $1 \leq s \leq n!$. 
The expected number of correct guesses of the deck after $k$ shuffles is
\begin{equation}
    e_k = E[X] = \sum_{s=1}^{n!} C^{(k)}_s \cdot P^{(k)}_{s},
\end{equation}  
where $P^{(k)}_{s}$ is the probability of getting the permutation $\pi_s$ after $k$ shuffles, and $C^{(k)}_s$ is the number of correct guesses of $\pi_s$ based on the optimal strategy. Here, 
$P^{(k)}_s = \dfrac{\binom{n+2^k-m_s}{n}}{2^{kn}}$, while the optimal strategy leading to $C^{(k)}_s$ has to be calculated specifically for each $n$ and $k$ following the program provided in the previous section. We give a simple example below.
}

{\color{black}
\textbf{Example 3} For $n=3$, we let $\pi_1 = [1,2,3], \pi_2=[1,3,2], \pi_3=[2,1,3], \pi_4=[2,3,1], \pi_5=[3,1,2]$ and $\pi_6=[3,2,1]$. Then $m_1$, the number of increasing subsequences of $\pi_1$, is 1 while $m_2=m_3=m_4=m_5= 2$ and $m_6=3$. Finally, it can be shown that the optimal strategy for a 3-card deck is the same for all $k$, i.e. guess the smallest number that has not shown up yet, and as a result $C^{(k)}_1=3, C^{(k)}_2=C^{(k)}_3=C^{(k)}_5=2, C^{(k)}_4=C^{(k)}_6=1$.
}

\textbf{Future work}

On one hand, finding an optimal guessing strategy for $k\geq2$ shuffles, or seeing if it even exists could be {\color{black}a crucial stepping stone} towards one direction for future work. 
On the other hand, we believe that tackling the card guessing version for $k \geq 2$, {\color{black}with no response (as considered by Ciucu in \cite{C})}, where the correctness of all the guesses is only revealed at the end of the game, seems to be a more promising direction for future work.

\subsection*{Acknowledgments}
The authors would like to thank Dr. Adam Chacon for constructive feedback of the manuscript.

\subsection*{Appendix A: Asymptotic approximation of $\binom{2n}{n}$\label{refA}}

This section of the Appendix is intended to be a handy reference to the standard method used to derive the infinite series formula for the central binomial. 

Starting with
\[ \ln n! \approx \ln\left( \sqrt{2\pi}n^{n+\frac{1}{2}}e^{-n} \right) 
+ \sum_{i=1}^{\infty}\frac{B_{2i}}{2i(2i-1)n^{2i-1}},\]
where $B_k$ is a Bernoulli number, we have
$$\ln{\binom{2n}{n}} = \ln(2n)! -2\ln n! 
= \ln\left( \frac{4^n}{\sqrt{\pi n}} \right) 
-\sum_{i=1}^{\infty}\frac{B_{2i}}{i(2i-1)n^{2i-1}} 
\left[ 1-1/2^{2i} \right] .$$

Finally, by Taylor series, we obtain
\[ \binom{2n}{n} \approx  \frac{4^n}{\sqrt{\pi n}} 
\left(1- \dfrac{1}{8n} +\dfrac{1}{128n^2}
+\dfrac{5}{1024n^3}-\dfrac{21}{32768n^4}+\dots   \right).  \]

\subsection*{Appendix B: Summary of Maple programs\label{refB}}
This section summarizes main programs used in each section. All Maple programs accompanying this paper are provided at \small{\url{https://thotsaporn.com/Card1.html}}.

\begin{small}
\textbf{Section 1: Introduction}

\textbf{Decks(n)}   \\
Input: Positive integer $n$ \\
Output: List of deck of $n$ cards after 1 shuffle \\ 
Try:  Decks(3);

\textbf{CorrectD(D)}   \\ 
Input: List of deck $D$ from one shuffle      \\
Output: Number of correct guesses of deck $D$ \\
Try: CorrectD([1,2,3,4]);

\end{small}

\begin{small}
\textbf{Section 2: One-time Riffle Shuffle}

\textbf{GenD(n,q)} \\ 
Input: Positive integer $n$ and symbolic $q$ \\
Output: Generating function of deck
of $n$ cards calculated by recurrence \eqref{D} \\
Try: GenD(15,q);

\textbf{CloseGr(r,k)} \\
Input: Positive integer $r$ and symbolic $k$ \\
Output: Closed form for $G^{(r)}_n(q)|_{q=1}$  in terms of $k$
for odd and even $n$ in Proposition \ref{prop:Gr} \\ 
Try: CloseGr(3,k);

\textbf{ForDr(r,s,K,n)} \\
Input: Positive integer $r$, parity 0 or 1 for $s$, positive integer $K \geq r$
  and symbolic $n$  \\
Output:   Formula in $n$ for $D^{(r)}_n(q)|_{q=1}$
for odd $n$ $(s=1)$ or even $n$ $(s=0)$ 
with about $K$ terms of accuracy \\  
Try: ForDr(1,1,5,n);

\textbf{CheckDr(n,r,K)} \\
Input: Positive integer $n, r$ and $K$ with $K \geq r$   \\
Output: Relative error of the value obtained with ForDr(r,s,K,n) to the exact value\\ 
Try: CheckDr(150,1,15);

\textbf{ForFacM(r,s,K,n)} \\
Input: Positive integer $r$, parity 0 or 1 for $s$, positive integer $K \geq r$
  and symbolic $n$  \\
Output:  Formula of factorial moment $E[(X)_r]$ in $n$
for odd $n$ $(s=1)$ or even $n$ $(s=0)$ 
with about $K$ terms of accuracy \\  
Try: ForFacM(3,1,5,n);

\textbf{MoMean(r,s,K,n)} \\
Input: Positive integer $r$, parity 0 or 1 for $s$, positive integer $K \geq r$
  and symbolic $n$  \\
Output:  Formula of moment about the mean $E[(X-\mu)^r]$ in $n$
for odd $n$ $(s=1)$ or even $n$ $(s=0)$ 
with about $K$ terms of accuracy \\  
Try: MoMean(3,1,5,n);
\end{small}

\begin{small}
\textbf{Section 3: Generalization to $k$ Riffle Shuffles: Is it possible?}

\textbf{KShuff(n,k)} \\
Input: Positive integers $n$ and $k$ \\
Output: List of deck of $n$ cards, with multiplicities, 
after shuffling $k$ times \\
Try: KShuff(4,4);

{\color{black}
\textbf{SlowBG(n,k,S), FastBG(n,k,S)} \\
Input: Positive integer $n$ and $k$, the list of all revealed cards $S$\\
Outputs: 1. The number that is most likely to show up next (after $S$) from the deck of $n$ cards with $k$ shuffles; 2. The list containing the numbers of ways that card number $i, \, 1 \leq i \leq n$ will be the next card (after $S$) after $k$ riffle shuffles.\\
Try: FastBG(5,2,[3]); or 
SlowBG(10,2,[5]); or FastBG(10,2,[5]); 
}
\end{small}

\begin{small}
\textbf{Section 4: Last But Not Least}

\textbf{Dora2(k),  Dora3(k)} \\
Input: Positive integer $k$ \\
Output: The expected number of correct
guesses after $k$ shuffles of deck of 2 (or 3) cards \\
Try: [seq(Dora3(k),k=0..10)];
\end{small}

\end{document}